\documentclass{amsart}
\usepackage{graphicx}
\usepackage{amsmath,amsthm,amssymb,IMjournal}
\newcommand{\p}{\partial}

\begin{document}
\newtheorem{The}{Theorem}[section]
\newtheorem{Lem}[The]{Lemma}
\newtheorem{Cor}[The]{Corollary}

\numberwithin{equation}{section}

\title{On the optimality and sharpness of Laguerre's lower bound on the smallest eigenvalue of a symmetric positive definite matrix}

\author{\|Yusaku |Yamamoto|, Chofu, Tokyo, 182-8585, Japan}

\rec {January 31, 2017}

\dedicatory{Cordially dedicated to ...}

\abstract 
   Lower bounds on the smallest eigenvalue of a symmetric positive definite matrices $A\in\mathbb{R}^{m\times m}$ play an important role in condition number estimation and in iterative methods for singular value computation. In particular, the bounds based on ${\rm Tr}(A^{-1})$ and ${\rm Tr}(A^{-2})$ attract attention recently because they can be computed in $O(m)$ work when $A$ is tridiagonal. In this paper, we focus on these bounds and investigate their properties in detail. First, we consider the problem of finding the optimal bound that can be computed solely from ${\rm Tr}(A^{-1})$ and ${\rm Tr}(A^{-2})$ and show that so called Laguerre's lower bound is the optimal one in terms of sharpness. Next, we study the gap between the Laguerre bound and the smallest eigenvalue. We characterize the situation in which the gap becomes largest in terms of the eigenvalue distribution of $A$ and show that the gap becomes smallest when ${\rm Tr}(A^{-2})/\{{\rm Tr}(A^{-1})\}^2$ approaches 1 or $\frac{1}{m}$. These results will be useful, for example,  in designing efficient shift strategies for singular value computation algorithms.
\endabstract



\thanks
   This study has been supported by JSPS KAKENHI Grant Numbers JP26286087, JP15H02708, JP15H02709 and JP16KT0016.
\endthanks

\section{Introduction}
\label{sec1}
Let $A \in{\bf R}^{m\times m}$ be a symmetric positive definite matrix and denote the smallest eigenvalue of $A$ by $\lambda_m(A)$. In this paper, we are interested in a lower bound on $\lambda_m(A)$. If the Cholesky factorization of $A$ is $A=BB^{\top}$, where $B\in{\bf R}^{m\times m}$ is a nonsingular lower triangular matrix, the smallest singular value of $B$ can be written as $\sigma_m(B)=\sqrt{\lambda_m(A)}$. Hence, finding a lower bound on $\lambda_m(A)$ is equivalent to finding a lower bound on $\sigma_m(B)$.

A lower bound on $\lambda_m(A)$ or $\sigma_m(B) $ plays an important role in various scientific computations. For example, when combined with an upper bound on $\|A\|_2$, a lower bound on $\lambda_m(A)$ can be used to give an upper bound on the condition number of $A$. In singular value computation algorithms such as the dqds algorithm \cite{Fernando94}, the orthogonal qd algorithm \cite {Matt97} and the mdLVs algorithm \cite{Iwasaki06}, a lower bound on $\sigma_m(B)$ is used as a shift to accelerate the convergence. In the latter case, the matrix $B$ is usually a lower bidiagonal matrix as a result of preprocessing by the Householder method \cite{Golub12}.

Several types of lower bounds on $\lambda_m(A)$ or $\sigma_m(B)$ have been proposed so far. There are bounds based on eigenvalue inclusion theorems such as Gershgorin's circle theorem \cite{Johnson89} or Brauer's oval of Cassini \cite{Johnson98}. The norm of the inverse, $\|A^{-1}\|_{\infty}$, can also be used to bound the maximum eigenvalue of $A^{-1}$ from above, and therefore to bound $\lambda_m(A)$ from below. There are also bounds based on the traces of the inverses, namely, ${\rm Tr}(A^{-1})$ and ${\rm Tr}(A^{-2})$. Among them, the last class of bounds are attractive in the context of singular value computation, because they always give a valid (positive) lower bound, as opposed to the bounds based on the eigenvalue inclusion theorems, and they can be computed in $O(m)$ work using efficient algorithms \cite{Kimura11, Yamashita12, Yamashita15}. Examples of lower bounds of this type include the Newton bound \cite {Matt97}, the generalized Newton bound \cite{ Kimura11,Aishima10} and the Laguerre bound \cite{Matt97}.

In this paper, we focus on the lower bounds of $\lambda_m(A)$ derived from ${\rm Tr}(A^{-1})$ and ${\rm Tr}(A^{-2})$ and investigate their properties. In particular, we will address the following two questions. The first is to identify an {\it optimal} formula for a lower bound on $\lambda_m(A)$ that is based solely on ${\rm Tr}(A^{-1})$ and ${\rm Tr}(A^{-2})$. Here, the word "optimal" means that the formula always gives a sharper (that is, larger) bound than any other formulas using only ${\rm Tr}(A^{-1})$ and ${\rm Tr}(A^{-2})$. As a result of our analysis, we show that the Laguerre bound mentioned above is the optimal formula in this sense. The second question is to evaluate the {\it gap} between the Laguerre bound and $\lambda_m(A)$. Unlike the Laguerre bound, $\lambda_m(A)$ is not determined uniquely only from ${\rm Tr}(A^{-1})$ and ${\rm Tr}(A^{-2})$. Hence, for some of the matrices, there must be a gap between the bound and $\lambda_m(A)$. Our problem is to quantify the maximum possible gap and identify the conditions under which the maximum gap is attained. These results will be useful, for example, in designing an efficient shift strategy for singular value computation algorithms, which combines the Laguerre bound with other bounds with complementary characteristics \cite{Yamashita13}.

The rest of this paper is structured as follows. In Section 2, we investigate the lower bounds on $\lambda_m(A)$ derived from ${\rm Tr}(A^{-1})$ and ${\rm Tr}(A^{-2})$ and show that the Laguerre bound is an optimal one in terms of sharpness. Section 3 deals with the gap between the Laguerre bound and $\lambda_m(A)$. In particular, we characterize the situation in which the gap becomes largest in terms of the eigenvalue distribution of $A$. Section 4 gives some concluding remarks.

\section{An optimal lower bound based on ${\rm Tr}(A^{-1})$ and ${\rm Tr}(A^{-2})$}
\subsection{Lower bounds based on ${\rm Tr}(A^{-1})$ and ${\rm Tr}(A^{-2})$}
Let $A$ be an $m\times m$ real symmetric positive matrix. We denote the $k$th largest eigenvalue of $A$ by $\lambda_k(A)$, or $\lambda_i$ for short. Let $f(\lambda)=\det(\lambda I-A)$ be the characteristic polynomial of $A$. To find a lower bound on the smallest eigenvalue $\lambda$, we consider applying a root finding method for an algebraic equation to $f(\lambda)=0$ starting from the initial value $\lambda^{(0)}=0$. There are several root finding methods, such as the Bailey's (Halley's) method \cite{Alefeld81}, Householder's method \cite{Householder70} and Laguerre's method \cite{Wilkinson88, Matt97}, for which the iteration formulas can be written as follows:
\begin{eqnarray}
\lambda_B^{(n+1)} &=& \lambda^{(n)}-\frac{f(\lambda^{(n)})}{f'(\lambda^{(n)})}\cdot\frac{1}{1-\frac{f(\lambda^{(n)})f''(\lambda^{(n)})}{2f'(\lambda^{(n)})^2}}
\label{eq:Bailey}, \\
\lambda_H^{({n+1})} &=& \lambda^{(n)}-\frac{f(\lambda^{(n)})}{f'(\lambda^{(n)})}\left\{1+\frac{f(\lambda^{(n)})f''(\lambda^{(n)})}{2f'(\lambda^{(n)})^2}\right\}
\label{eq:Householder}, \\
\lambda_L^{({n+1})} &=& \lambda^{(n)}-\frac{f(\lambda^{(n)})}{f'(\lambda^{(n)})} \nonumber \\
&& \quad \quad\times\frac{m}{1+\sqrt{(m-1)\left\{m\cdot\frac{f'(\lambda^{(n)})^2-f(\lambda^{(n)})f''(\lambda^{(n)})}{f'(\lambda^{(n)})^2}-1\right\}}}
\label{eq:Laguerre}
\end{eqnarray}
Eqs.~(\ref{eq:Bailey}), (\ref{eq:Householder}) and (\ref{eq:Laguerre}) represent the iteration formulas of Bailey's method, Householder's method and Laguerre's method, respectively. When applied to $f(\lambda)=\det(\lambda I-A)$ starting from $\lambda^{(0)}=0$, these formulas produce a sequence that increases monotonically and converges to $\lambda_m$. Hence, all of $\lambda_B^{(1)}$, $\lambda_H^{(1)}$ and $\lambda_L^{(1)}$ can be used as a lower bound on $\lambda_m$.

Noting that $f(\lambda)=\prod_{k=1}^m(\lambda-\lambda_k)$, we have
\begin{eqnarray}
f'(\lambda) &=& -\sum_{k=1}^m \prod_{j\ne k}(\lambda_j-\lambda) \nonumber \\
&=& -\prod_{j=1}^m(\lambda_j-\lambda)\sum_{k=1}^m\frac{1}{\lambda_k-\lambda} = -f(\lambda){\rm Tr}\left(\left(A-\lambda I\right)^{-1}\right), \\
f''(\lambda) &=& -f'(\lambda){\rm Tr}\left(\left(A-\lambda I\right)^{-1}\right)-f(\lambda)\sum_{k=1}^m\frac{1}{(\lambda_k-\lambda)^2} \nonumber \\
&=& -f'(\lambda){\rm Tr}\left(\left(A-\lambda I\right)^{-1}\right)-f(\lambda){\rm Tr}\left(\left(A-\lambda I\right)^{-2}\right).
\end{eqnarray}
Hence,
\begin{eqnarray}
\frac{f(\lambda)}{f'(\lambda)} &=& -\frac{1}{{\rm Tr}\left(\left(A-\lambda I\right)^{-1}\right)}, \\
\frac{f(\lambda)f''(\lambda)}{f'(\lambda)^2} &=& 1-\frac{{\rm Tr}\left(\left(A-\lambda I\right)^{-2}\right)}{\left\{{\rm Tr}\left(\left(A-\lambda I\right)^{-1}\right)\right\}^2}.
\end{eqnarray}
Inserting these into Eqs.~(\ref{eq:Bailey}), (\ref{eq:Householder}) and (\ref{eq:Laguerre}) with $\lambda^{(0)}=0$, we obtain the following lower bounds on $\lambda_m(A)$:
\begin{eqnarray}
L_B(A) &=& \frac{2 {\rm Tr}\left(A^{-1}\right)}{\left\{{\rm Tr}\left(A^{-1}\right)\right\}^2+{\rm Tr}\left(A^{-2}\right)}, \\
L_H(A) &=& \frac{1}{{\rm Tr}\left(A^{-1}\right)}\left[\frac{3}{2}-\frac{1}{2}\cdot\frac{{\rm Tr}\left(A^{-2}\right)}{\left\{{\rm Tr}\left(A^{-1}\right)\right\}^2}\right], \\
L_L(A) &=& \frac{1}{{\rm Tr}\left(A^{-1}\right)}\cdot\frac{m}{1+\sqrt{(m-1)\left[m\cdot\frac{{\rm Tr}\left(A^{-2}\right)}{\left\{{\rm Tr}\left(A^{-1}\right)\right\}^2}-1\right]}}. \label{eq:Laguerre3}
\end{eqnarray}
We call $L_B(A)$, $L_H(A)$ and $L_L(A)$ the Bailey bound, the Householder bound and the Laguerre bound, respectively. In addition to these, we also have a simple bound:
\begin{equation}
L_N(A) = \{{\rm Tr}(A^{-2})\}^{-\frac{1}{2}} \le \left(\sum_{k=1}^m\frac{1}{\lambda_k^2}\right)^{-\frac{1}{2}} < \lambda_m,
\end{equation}
which is called the Newton bound of order 2 \cite{Matt97, Kimura11, Aishima10}. In the case where $A$ is a tridiagonal matrix, both ${\rm Tr}(A^{-1})$ and ${\rm Tr}(A^{-2})$ can be computed in $O(m)$ work from its Cholesky factor $B$ \cite{Kimura11, Yamashita12, Yamashita15}. Accordingly, any of these bounds can be employed in a practical shift strategy for singular value computation algorithms. The problem then is which of the four lower bounds, or possibly another bound derived from ${\rm Tr}(A^{-1})$ and ${\rm Tr}(A^{-2})$, is optimal in terms of sharpness.

\subsection{The optimal lower bound}
To answer the question, we reformulate the problem as follows. Assume that ${\rm Tr}(A^{-1})$ and ${\rm Tr}(A^{-2})$ are specified for a symmetric positive definite matrix $A$. Then, how small can the smallest eigenvalue $\lambda_m(A)$ be? If this bound can be obtained explicitly as a function of ${\rm Tr}(A^{-1})$ and ${\rm Tr}(A^{-2})$, then, it will be the optimal formula for the lower bound of $\lambda_m(A)$.

Now, let $a\equiv{\rm Tr}(A^{-1})$, $b\equiv{\rm Tr}(A^{-2})$ and $x_k\equiv 1/\lambda_k$ ($k=1, 2, \ldots, m$). Then, the upper bound on $x_m$ (the reciprocal of the lower bound on $\lambda_m$) can be obtained by solving the following constrained optimization problem:
\begin{eqnarray}
{\rm maximize} \;\; x_m && \label{eq:objective} \\
{\rm s.t.} && \sum_{k=1}^m x_k = a, \label{eq:constraint1} \\
&& \sum_{k=1}^m x_k^2 = b, \label{eq:constraint2} \\
&& x_k > 0 \quad(k=1, 2, \ldots, m), \label{eq:constraint3} \\
&& x_1 \le x_2 \le \cdots \le x_m. \label{eq:constraint4}
\label{eq:inequality}
\end{eqnarray}
Actually, the constraint (\ref{eq:constraint4}) is redundant, because if $(x_1^*, x_2^*, \ldots, x_m^*)$ is a solution of the optimization problem without constraint (\ref{eq:constraint4}), then, from symmetry, $(x_{\sigma(1)}^*, x_{\sigma(2)}^*, \ldots, x_{\sigma(m)}^*)$ is also a solution for any permutation $\sigma$ of $\{1, 2, \ldots, m\}$, and therefore we can choose a solution that satisfies (\ref{eq:constraint4}). Hence we omit (\ref{eq:constraint4}) in the following.

To solve the optimization problem (\ref{eq:objective})--(\ref{eq:constraint3}), we remove the constraint (\ref{eq:constraint3}) and consider a relaxed problem described by (\ref{eq:objective})--(\ref{eq:constraint2}). By introducing the Lagrange multipliers $\mu$ and $\nu$, we can write the Lagrangian as
\begin{equation}
L = x_m - \mu\left(\sum_{k=1}^m x_k - a\right) - \nu\left(\sum_{k=1}^m x_k^2 - b\right).
\end{equation}
Then the solution to (\ref{eq:objective})--(\ref{eq:constraint2}) must satisfy
\begin{eqnarray}
\frac{\p L}{\p x_m} &=& 1-\mu-2\nu x_m = 0, \label{eq:L1} \\
\frac{\p L}{\p x_k} &=& -\mu-2\nu x_k = 0 \quad(k=1, 2, \ldots, m-1), \label{eq:L2} \\
\frac{\p L}{\p \lambda} &=& \sum_{k=1}^m x_k - a = 0, \label{eq:L3} \\
\frac{\p L}{\p \mu} &=& \sum_{k=1}^m x_k^2 - b = 0. \label{eq:L4}
\end{eqnarray}
From (\ref{eq:L2}), we have either $\nu=0$ or $x_1=x_2=\cdots=x_{m-1}$. However, when $\nu=0$, we have $\mu=0$ from (\ref{eq:L2}) and $\mu=1$ from (\ref{eq:L1}), which is a contradiction. Thus $x_1=x_2=\cdots=x_{m-1}$ must hold. Inserting this into (\ref{eq:L3}) and (\ref{eq:L4}) leads to
\begin{eqnarray}
&& x_m + (m-1)x_1 - a = 0, \label{eq:xm1} \\
&& x_m^2 + (m-1)x_1^2 - b = 0.
\end{eqnarray}
Solving these simultaneous equations with respect to $x_m$ gives
\begin{equation}
x_m^{\pm} = \frac{a\pm\sqrt{m(m-1)b-(m-1)a^2}}{m} \label{eq:optimal}.
\end{equation}
Note that the $x_m$ given by (\ref{eq:optimal}) is real, since
\begin{eqnarray}
m(m-1)b-(m-1)a^2 &=& (m-1)\left\{m\sum_{k=1}^m x_k^2-\left(\sum_{k=1}^m x_k\right)^2\right\} \nonumber \\
&=& (m-1)\sum_{k=1}^m\sum_{l=1}^{k-1}(x_k-x_l)^2 \ge 0.
\label{eq:alphalower}
\end{eqnarray}

Now we return to the relaxed optimization problem (\ref{eq:objective})--(\ref{eq:constraint2}). Since the feasible set of this problem is compact and both the objective function and the constraints are differentiable, it must have a minimum and a maximum at a point where the gradient of the Lagrangian is zero. Furthermore, since the objective function is $x_m$ itself, the maximum is attained when $x_m=x_m^+$. Then, from Eq.~(\ref{eq:xm1}), we have
\begin{equation}
x_1 = x_2 = \cdots = x_{m-1} = \frac{(m-1)a-\sqrt{m(m-1)b-(m-1)a^2}}{m(m-1)}. \label{eq:xkoptimal}
\end{equation}
Hence, Eq.~(\ref{eq:xkoptimal}) and $x_m=x_m^+$ are the solution of the relaxed optimization problem.

Finally, we consider the positivity constraint (\ref{eq:constraint3}). It is clear from (\ref{eq:optimal}) that $x_m^+>0$. To investigate the positivity of the other variables, note that
\begin{equation}
a^2-b = \left(\sum_{k=1}^m\frac{1}{\lambda_k}\right)^2-\sum_{k=1}^m\frac{1}{\lambda_k^2} = 2\sum_{k=1}^m \sum_{l=1}^{k-1} \frac{1}{\lambda_k}\cdot\frac{1}{\lambda_l}>0,
\label{eq:alphaupper}
\end{equation}
where we used the fact that $a$ and $b$ are the traces of the inverse of a matrix with positive eigenvalues. Then (\ref{eq:xkoptimal}) can be rewritten as
\begin{eqnarray}
x_1 = x_2 = \cdots = x_{m-1} 
&=& \frac{(m-1)^2a^2-\{m(m-1)b-(m-1)a^2\}}{m(m-1)\{(m-1)a+\sqrt{m(m-1)b-(m-1)a^2}\}} \nonumber \\
&=& \frac{m(m-1)(a^2-b)}{m(m-1)\{ (m-1)a+\sqrt{m(m-1)b-(m-1)a^2}\}}>0.
\end{eqnarray}
This shows that the solution to the relaxed problem (\ref{eq:objective})--(\ref{eq:constraint2}) automatically satisfies the constraint (\ref{eq:constraint3}). Hence it is also a solution to the original problem (\ref{eq:objective})--(\ref{eq:constraint3}). Returning to the original variables $\lambda_k=1/x_k$, we know that the smallest value that $\lambda_m$ can take is
\begin{equation}
\frac{1}{{\rm Tr}\left(A^{-1}\right)}\cdot\frac{m}{1+\sqrt{(m-1)\left[m\cdot\frac{{\rm Tr}\left(A^{-2}\right)}{\left\{{\rm Tr}\left(A^{-1}\right)\right\}^2}-1\right]}}.
\label{eq:Laguerre2}
\end{equation}
This gives the optimal lower bound on $\lambda_m(A)$ in terms of ${\rm Tr}(A^{-1})$ and ${\rm Tr}(A^{-2})$. Since Eq.~(\ref{eq:Laguerre2}) is exactly the Laguerre bound (\ref{eq:Laguerre3}), we arrive at the following theorem.
\begin{The}
\label{T1}
Among the lower bounds on $\lambda_m(A)$ computed from ${\rm Tr}(A^{-1})$ and ${\rm Tr}(A^{-2})$, the Laguerre bound (\ref{eq:Laguerre3}) is optimal in terms of sharpness.
\end{The}

\section{The gap between the Laguerre bound and the smallest eigenvalue}
Now that we have established that the Laguere bound is the optimal lower bound, we next study the gap between the bound and the minimum eigenvalue. We begin with a lemma that holds for a $3\times 3$ matrix and then proceed to the general case. In the course of discussion, we also allow infinite eigenvalues to make the arguments simpler.

Assume that $A\in{\mathbb R}^{3\times 3}$ is a symmetric positive definite matrix with ${\rm Tr}(A^{-1})=a$ and ${\rm Tr}(A^{-2})=b$. Let the eigenvalues of $A$ be $\lambda_1\ge \lambda_2\ge \lambda_3>0$. To evaluate the gap, we consider how large $\lambda_3$ can be under the fixed values of $a$ and $b$. First, we show the following lemma.
\begin{Lem}
\label{Lem2}
For fixed $a={\rm Tr}(A^{-1}) $ and $b={\rm Tr}(A^{-2})$, $\lambda_3$ can take a maximum only when $\lambda_2=\lambda_3$ or $\lambda_1=+\infty$.
\end{Lem}
\proof
Let $x=1/\lambda_3$, $y=1/\lambda_2$ and $z=1/\lambda_1$. Since we allow infinite eigenvalues, the point $(x, y, z)$ lies in a region $D$ of the $xyz$ space specified by $x+y+z=a$, $x^2+y^2+z^2=b$, and $x\ge y\ge z\ge 0$. Since $D$ is a compact set, the continuous function $x$ attains a minimum somewhere in $D$. Hence, if we can show that $x$ does not attain a minimum when $x>y$ and $z>0$, it means that $x$ attains a minimum when $x=y$ or $z=0$.

Assume that the point $(x,y,z)$ is in $D$ and both $x>y$ and $z>0$ hold. Then, let $\epsilon>0$ be some small quantity and $t\in{\bf R}$ and consider changing $(x,y,z)$ to $(x^{\prime}, y^{\prime}, z^{\prime})$ as follows:
\begin{eqnarray}
x^{\prime} &=& x - \epsilon, \label{eq:xprime} \\
y^{\prime} &=& y + t\epsilon, \label{eq:yprime} \\
z^{\prime} &=& z + (1-t)\epsilon. \label{eq:zprime} 
\end{eqnarray}
Clearly, the new point $(x^{\prime}, y^{\prime}, z^{\prime})$ lies on the plane $x+y+z=a$. We determine $t$ so that it is also on the sphere $x^2+y^2+z^2=b$. The condition can be written as
\begin{equation}
(x-\epsilon)^2 + (y+t\epsilon)^2 + \{z+(1-t)\epsilon\}^2 = x^2+y^2+z^2,
\end{equation}
or
\begin{equation}
\epsilon t^2 + (y-z-\epsilon)t + (-x+z+\epsilon)=0.
\end{equation}
Solving this with respect to $t$ gives
\begin{equation}
t_{\pm}= \frac{-(y-z-\epsilon)\pm\sqrt{(y-z-\epsilon)^2+4\epsilon(x-z-\epsilon)}}{2\epsilon}.
\label{eq:tsolution}
\end{equation}
In the following, we adopt the solution $t=t_+$. Now we consider two cases. First, consider the case of $y=z$. Then, we have from Eq.~(\ref{eq:tsolution}),
\begin{equation}
t_+\epsilon = \frac{\epsilon+\sqrt{\epsilon^2+4\epsilon(x-z-\epsilon)}}{2} = O(\sqrt{\epsilon}).
\end{equation}
Inserting this into (\ref{eq:xprime}) through (\ref{eq:zprime}), we know that the changes in $x$, $y$ and $z$ are at most $O(\sqrt{\epsilon})$ when $\epsilon$ is small.

Next, consider the case of $y>z$. In this case, we can rewrite (\ref{eq:tsolution}) as
\begin{equation}
t_+ = \frac{2(x-z-\epsilon)}{(y-z-\epsilon)+\sqrt{(y-z-\epsilon)^2+4\epsilon(x-z-\epsilon)}}.
\end{equation}
Since $x-z-\epsilon>0$ and $y-z-\epsilon>0$ for a sufficiently small $\epsilon$, we have
\begin{equation}
0 < t_+ < \frac{x-z-\epsilon}{y-z-\epsilon}.
\end{equation}
Hence, $ t_+\epsilon=O(\epsilon)$ when $\epsilon$ is small and therefore the changes in $x$, $y$ and $z$ are at most $O(\epsilon)$ in this case.

\begin{figure}[h]
\centerline{\includegraphics[height=1.0in]{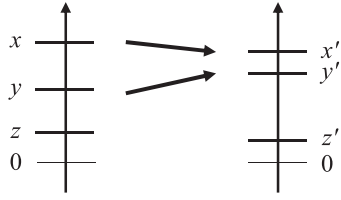}}
\caption{The values of $x$, $y$ and $z$ before and after the perturbation.}
\label{fig:fig1}
\end{figure}

In summary, in both cases, the changes of $x$, $y$ and $z$ can be made arbitrarily small. Thus, by choosing $\epsilon$ sufficiently small, we can make $x^{\prime}$ smaller than $x$ while keeping the relation $x^{\prime}>y^{\prime}>0$ and $x^{\prime}>z^{\prime}>0$ (Fig.~\ref{fig:fig1}). The relation $y^{\prime}\ge z^{\prime}$ may not hold, but in that case, we can interchange $y^{\prime}$ and $z^{\prime}$. In this way, we can obtain another point $(x^{\prime}, y^{\prime}, z^{\prime})\in D$ which attains a smaller value of $x$. Hence $x$ cannot attain a minimum when both $x>y$ and $z>0$ hold and the lemma is proved.
\endproof

Using this lemma, we can prove the following theorem.
\begin{The}
\label{T3}
Let $a={\rm Tr}(A^{-1}) $ and $b={\rm Tr}(A^{-2})$ be fixed and $q$ be an integer satisfying $q<a^2/b\le q+1$. Then, $\lambda_m(A)$ takes a maximum when $\lambda_1(A)=\cdots=\lambda_{m-q-1}(A)=+\infty$ and $\lambda_{m-q+1}=\cdots=\lambda_m(A)$. The maximum is given as
\begin{equation}
\lambda_m^*(A) = \frac{1}{{\rm Tr}\left(A^{-1}\right)}\cdot\frac{q(q+1)}{q+\sqrt{q\left[(q+1)\cdot\frac{{\rm Tr}\left(A^{-2}\right)}{\left\{{\rm Tr}\left(A^{-1}\right)\right\}^2}-1\right]}}.
\label{eq:Theorem3}
\end{equation}

\end{The}
\proof
Let $x_k=1/\lambda_k$. First, assume that there are two or more eigenvalues which are neither an infinite eigenvalue nor equal to $\lambda_m(A)$. In this case, as we will show in the following, we can make $\lambda_m$ smaller by adding appropriate perturbations. We divide the cases depending on the multiplicity $q$ of the smallest eigenvalue.

\begin{figure}[h]
\centerline{\includegraphics[height=1.1in]{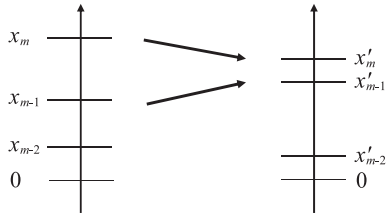}}
\caption{The values of $x_m$, $x_{m-1}$ and $x_{m-2}$ before and after the perturbation.}
\label{fig:fig2}
\end{figure}

When $q=1$, from the assumption, both $\lambda_{m-2}$ and $\lambda_{m-1}$ are neither an infinite eigenvalue nor equal to $\lambda_m(A)$. Thus, we have $0<x_{m-2}\le x_{m-1}<x_m$. Then, by picking up these three variables and adding the same perturbations as in Lemma \ref{Lem2}, we can make $x_m$ smaller while keeping the condition $0<x_{m-2}\le x_{m-1}<x_m$ (Fig.~\ref{fig:fig2}). Clearly, the values of $ {\rm Tr}(A^{-1}) $ and $ {\rm Tr}(A^{-2})$ are unchanged by this perturbation. Hence, $x_m$ cannot take a minimum in this case.

\begin{figure}[h]
\centerline{\includegraphics[height=1.1in]{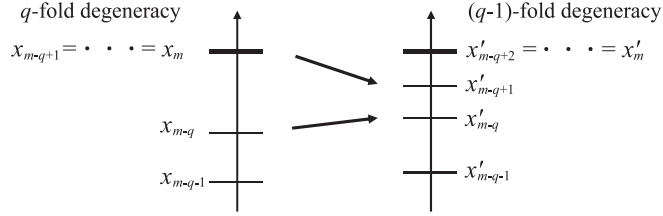}}
\caption{The values of $x_{m-q-1}, x_{m-q}, \ldots, x_m$ before and after the perturbation.}
\label{fig:fig3}
\end{figure}

When $q>1$, $0<x_{m-q-1}\le x_{m-q}<x_{m-q+1}=\cdots=x_m$ holds from the assumption. Then, by picking up the three variables $x_{m-q-1}$, $x_{m-q}$ and $x_{m-q+1}$ and adding the perturbations as in Lemma \ref{Lem2}, we can make $x_{m-q+1}$ smaller while keeping $0<x_{m-q-1}\le x_{m-q}<x_{m-q+1}$. This does not change the smallest eigenvalue, but reduces its multiplicity from $q$ to $q-1$ (Fig.~\ref{fig:fig3}). Moreover, the condition that there are two or more eigenvalues which are neither an infinite eigenvalue nor equal to $\lambda_m(A)$ still holds. Hence, we can repeat this procedure and reduce $q$ to 1, while keeping the value of the smallest eigenvalue unchanged. But in this last situation, $x_m$ cannot take a minimum, as concluded in the analysis of the $q=1$ case. 

From the above analysis, we can conclude that $x_m$ cannot take a minimum when there are two or more eigenvalues which are neither an infinite eigenvalue nor equal to $\lambda_m(A)$. Thus, the only possible case is when $x_1=\cdots=x_{m-q-1}=0$ and $x_{m-q+1}=\cdots=x_m$ holds for some $q$. In this case, we have
\begin{eqnarray}
x_{m-q}+q x_m &=& a, \\
x_{m-q}^2 + q x_m^2 &=& b,
\end{eqnarray}
or
\begin{eqnarray}
x_m^{\pm} &=& \frac{aq\pm\sqrt{q\{(q+1)b-a^2\}}}{q(q+1)}, \\
x_{m-q}^{\pm} &=& \frac{a\mp\sqrt{q\{(q+1)b-a^2\}}}{q+1}.
\end{eqnarray}
For $x_m$ and $x_{m-q}$ to be real, $q$ must satisfy $q+1\ge a^2/b$. Then, for $x_{m-q}\le x_m$ to hold, we have to choose $x_m^+$ and $x_{m-q}^+$. In addition, for $x_{m-q}^+>0$ to hold, we must have $q<a^2/b$. From the condition $q<a^2/b\le q+1$, $q$ is determined uniquely. Hence, there is only one set of $q$, $x_m$ and $x_{m-q}$ that satisfy the condition for minimum $x_m$. Since the feasible region of $(x_1, x_2, \ldots, x_m)$, specified by $\sum_{k=1}^m x_k=a$, $\sum_{k=1}^m x_k^2=b$ and $0\le x_1\le x_2\le \cdots \le x_m$, is compact, $x_m$ must have a minimum somewhere in this region. Accordingly, we conclude that $x_m$ takes a minimum when $q<a^2/b\le q+1$, $x_1=\cdots=x_{m-q-1}=0$, $x_{m-q}=x_{m-q}^+$ and $x_{m-q+1}=\cdots=x_m=x_m^+$. Eq.~(\ref {eq:Theorem3}) is obtained from $\lambda_m=1/x_m$.
\endproof

To measure the gap between the Laguerre bound and the smallest eigenvalue, we use the quantity $L_L(A)/\lambda_m^*(A)$, which becomes one when there is no gap and zero when the gap is maximal. Let $\alpha\equiv {\rm Tr}(A^{-2})/\{{\rm Tr}(A^{-1})\}^2$ and $q$ be an integer specified in Theorem \ref{T3}. Then, we have from Eqs.~(\ref{eq:Theorem3}) and (\ref{eq:Laguerre3}),
\begin{equation}
\frac{L_L(A)}{\lambda_m^*(A)} = \frac{m}{q(q+1)}\cdot\frac{q+\sqrt{q\{(q+1)\alpha-1\}}}{1+\sqrt{(m-1)(m\alpha-1)}}.
\label{eq:gap}
\end{equation}
Thus, we have obtained an expression for the maximum possible gap as a function of $m$ and $\alpha$ (note that $q$ is determined from $\alpha$ uniquely).

So far, we have allowed infinite eigenvalues. However, of course, actual matrices have only finite eigenvalues. Accordingly, except for the case of $q=m-1$, for which no infinite eigenvalues are required for $\lambda_m(A)$ to take a maximum, the right-hand side of (\ref{eq:gap}) is a lower bound that can be approached arbitrarily closely.

Finally, we investigate the behavior of the right-hand side of (\ref{eq:gap}) as a function of $\alpha$. Note that $\frac{1}{m}\le \alpha <1$ from (\ref{eq:alphalower}) and (\ref{eq:alphaupper}). We consider three extreme cases, namely, $\frac{1}{m}\le\alpha<\frac{1}{m-1}$, $\frac{1}{m}\ll\alpha\ll 1$ and $\frac{1}{2}\le\alpha<1$.
\begin{itemize}
\item When $\frac{1}{m}\le\alpha<\frac{1}{m-1}$, we have $q=m-1$ and therefore
\begin{equation}
\frac{L_L(A)}{\lambda_m^*(A)} =\frac{1}{m-1}\cdot\left\{1+\frac{m-2}{1+\sqrt{(m-1)(m\alpha-1)}}\right\}.
\end{equation}
This is a decreasing function in $\alpha$ and takes the maximum value 1 at $\alpha=\frac{1}{m}$ and approaches $\frac{1}{2}\cdot\frac{m}{m-1}$ as $\alpha\rightarrow\frac{1}{m-1}$. Hence $L_L(A)/\lambda_m^*(A)>\frac{1}{2}$ all over the region.
\item When $\frac{1}{m}\ll\alpha\ll 1$, we have $1\ll q\ll m$ and therefore $L_L(A)/\lambda_m^*(A)\simeq 1/\sqrt{q}$.
\item When $\frac{1}{2}\le\alpha<1$, we have $q=1$ and therefore
\begin{equation}
\frac{L_L(A)}{\lambda_m^*(A)} = \frac{m}{2}\cdot\frac{1+\sqrt{2\alpha-1}}{1+\sqrt{(m-1)(m\alpha-1)}}.
\end{equation}
For $\frac{1}{2}\le\alpha<1$, this is an increasing function in $\alpha$ that takes the minimum value
\begin{equation}
\frac{1}{\sqrt{2}}\cdot\frac{1}{\frac{\sqrt{2}}{m}+\sqrt{\left(1-\frac{1}{m}\right)\left(1-\frac{2}{m}\right)}}
\end{equation}
at $\alpha=\frac{1}{2}$ and approaches 1 as $\alpha\rightarrow 1$. Thus, when $m$ is large, $L_L(A)/\lambda_m^*(A)$ is almost larger than $\frac{1}{\sqrt{2}}$ all over the region.
\end{itemize}

In summary, we can conclude that the Laguerre bound is fairly tight when $\alpha$ is close to $\frac{1}{m}$ or greater than $\frac{1}{2}$ and can be loose when $\alpha$ is in the intermediate region.

\begin{figure}[h]
\centerline{\includegraphics[height=2.6in]{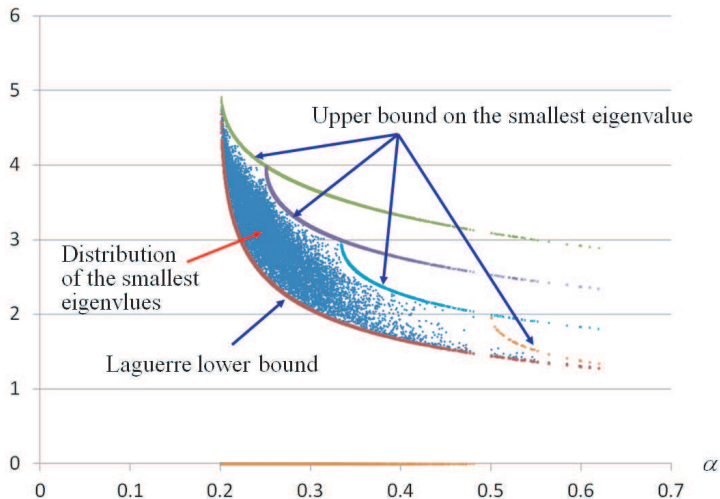}}
\caption{The smallest eigenvalues of randomly generated $5\times 5$ symmetric positive definite matrices as a function of $\alpha$.}
\label{fig:fig4}
\end{figure}

In Fig.~\ref{fig:fig4}, we plot the smallest eigenvalues of randomly generated $5\times 5$ symmetric positive definite matrices. These matrices are normalized so that ${\rm Tr}(A^{-1})=1$ and the horizontal axis is $\alpha={\rm Tr}(A^{-2})/\{{\rm Tr}(A^{-1})\}^2={\rm Tr}(A^{-2})$. The Laguerre bound (\ref{eq:Laguerre3}) and the upper bound (\ref{eq:Theorem3}) on the smallest eigenvalue are also shown in the graph. From the graph, we can confirm the optimality of the Laguerre bound, since it actually constitutes the lower boundary of the region where the smallest eigenvalues exist. We also see that the upper boundary is given by (\ref{eq:Theorem3}). Finally, it is clear that the Laguerre bound is tight when $\alpha\simeq\frac{1}{m}$ or $\alpha\ge\frac{1}{2}$ and loose in the intermediate region.

\section{Conclusion}
In this paper, we investigated the properties of lower bounds on the smallest eigenvalue of a symmetric positive definite matrix $A$ computed from ${\rm Tr}(A^{-1})$ and ${\rm Tr}(A^{-2})$. We studied two problems, namely, finding the optimal bound and evaluating its sharpness. As for the first question, we found that the Laguerre bound is the optimal one in terms of sharpness. As for the second question, We characterized the situation in which the gap becomes largest in terms of the eigenvalue distribution of $A$. Furthermore, we showed that the gap becomes smallest when ${\rm Tr}(A^{-2})/\{{\rm Tr}(A^{-1})\}^2$ approaches 1 or $\frac{1}{m}$. These results will help designing efficient shift strategies for singular value computation methods such as the dqds algorithm and the mdLVs algorithm.

\paragraph{\bf Acknowledgment}
The author would like to thank Prof.~Kinji Kimura, Dr.~Yuji Nakatsukasa and Dr.~Takumi Yamashita for fruitful discussion.

{\small
{\em Authors' addresses}:
{\em Yusaku Yamamoto}, The University of Electro-Communications, 1-5-1 Chofugaoka, Chofu, Tokyo, 182-8585, Japan.
 e-mail: \texttt{yusaku.yamamoto@uec.ac.jp}.

}

\end{document}